\documentclass[12pt]{amsart}

\pdfoutput=1

\usepackage{amssymb}
\usepackage{graphicx}

\advance\evensidemargin-.5in
\advance\oddsidemargin-.5in
\advance\textwidth1in

\makeatletter
\def\ps@firstpage{\ps@plain
  \def\@oddfoot{\normalfont\scriptsize \hfil\thepage\hfil
     \global\topskip\normaltopskip}%
  \let\@evenfoot\@oddfoot
  \def\@oddhead{\vbox{
\vskip-18pt
{\hbox to\hsize{
\hfill
arXiv:1002.4921 [math.AG]
}}
\vskip6pt
{\hbox to\hsize{
\hfill
UCSB Math 2010-02
}}
}}
  \let\@evenhead\@oddhead 
}
\makeatother

\newtheorem{conjecture}{Conjecture}

\theoremstyle{definition}

\theoremstyle{remark}

\numberwithin{equation}{section}

\begin{document}

\title[On the structure of supersymmetric $T^3$ fibrations]{On the structure 
of supersymmetric\\ $T^3$ fibrations}

\author{David R. Morrison}
\address{University of California, Santa Barbara}
\curraddr{}
\email{drm@math.ucsb.edu}
\thanks{This research was
supported in part by the  National Science Foundation (grant
DMS-0606578).
Any opinions, findings, and conclusions or recommendations expressed in this
material are those of the author
and do not necessarily reflect the views of the National Science Foundation.}

\date{}

\begin{abstract}
We formulate some precise conjectures concerning the
existence and structure of
supersymmetric $T^3$ fibrations of Calabi--Yau threefolds, and 
describe how these conjectural
fibrations would give rise to the Strominger--Yau--Zaslow version of
mirror symmetry.
\end{abstract}

\maketitle

Mirror symmetry between Calabi--Yau manifolds 
remains, some twen\-ty years after its discovery, one of the biggest mysteries
in mathematics.  Originally formulated as a physical relationship between
certain pairs of Calabi--Yau manifolds \cite{CLS,GP,Aspinwall:1990xe}
(with astonishing mathematical consequences
relating enumerative geometry to Hodge theory \cite{cdgp,mirrorguide}), 
the mirror symmetry proposal was
refined in 1996 by Strominger, Yau, and Zaslow \cite{syz}
into a much more geometric
statement, again based on physics.  But while the 
main idea of the Strominger--Yau--Zaslow
 proposal has been clear from the outset, many of
the details have remained elusive.  One of the reasons for this is that
the proposal involves special Lagrangian submanifolds of a Calabi--Yau
manifold, and very few tools are available for studying such submanifolds.
As a consequence, the initial 
period of intense study 
of the original Strominger--Yau--Zaslow
proposal has largely ended,\footnote{A very readable
summary of the progress made on the original proposal,
and the transition to the more recent approaches, was given
recently by Gross \cite{gross-SYZ}.}
and much of the recent work on mirror symmetry has
shifted to other approaches 
(including reformulations of the Strominger--Yau--Zaslow proposal),
as is recounted in detail elsewhere in this
volume.

Our purpose in this paper is to give a quite precise conjectural formulation
of the Strominger--Yau--Zaslow 
version of mirror symmetry for Calabi--Yau threefolds.  
The conjectures we formulate are modifications of conjectures previously
made by Gross, Ruan, and Joyce; we also restate some other conjectures from
\cite{syz,GrossWilson:LCSL,KS:mirror}.
Our
formulation is unfortunately not directly
based on examples, since---as mentioned
above---tools for constructing concrete examples are not currently available.
The conjectures are, however, motivated on the one hand by qualitative
features of the Strominger--Yau--Zaslow 
proposal which have been discovered by mathematicians,
and on the other hand by some suggestive
arguments from physics.
We focus on the motivation from mathematics in this paper.

For simplicity, 
we restrict our attention for the most part
to Calabi--Yau manifolds of
complex dimension one, two, or three.  We expect, however, that similar
conjectures could be formulated in higher dimension, at the expense of
greater combinatorial complexity.

\section{Supersymmetric torus fibrations}

A {\em Calabi--Yau metric}\/ 
is a Riemannian metric on a manifold $X$ of dimension
$2n$ whose Riemannian holonomy is precisely $SU(n)$.\footnote{Some authors
use the term {\em Calabi--Yau metric}\/  when the holonomy is any subgroup of
$SU(n)$.}  The representation theory of the holonomy group gives
rise to various geometric structures on $X$: 
there is a compatible complex structure
(unique up to complex conjugation when $n\ge3$), a $2$-form $\omega$ which
serves as the K\"ahler form
of the given
metric with respect to that complex structure, and a nowhere-vanishing
holomorphic $n$-form $\Omega$.  Both $\omega$ and $\Omega$ are covariantly
constant (as is the almost-complex structure operator $\mathcal J$).

Thanks to a result conjectured by Calabi \cite{Calabi} and proven by
Yau \cite{Yau}, when $X$ is a compact K\"ahler manifold with a
nowhere-vanishing holomorphic $n$-form,  there is a unique
Ricci-flat metric in each de Rham cohomology class 
$[\omega]\in H^2(X,\mathbb R)$ containing a K\"ahler form.
These
metrics have holonomy contained in $SU(n)$, so they will be Calabi--Yau
under our definition
provided that the holonomy does not reduce to a subgroup.  As a
consequence, the existence of a nowhere-vanishing holomorphic $n$-form
on a K\"ahler manifold $X$
is often taken as a definition of {\em Calabi--Yau manifold}.

Given a Calabi--Yau metric on $X$, an associated K\"ahler form $\omega$,
and a holomorphic $n$-form $\Omega$, we say that a submanifold $L$ of (real)
dimension $n$ is
{\em special Lagrangian}\/ if $\omega|_L\equiv0$ and 
$\operatorname{Im}(e^{i\theta}\Omega)|_L\equiv0$ for some $\theta$ called
the {\em phase}\/ of $L$. This notion was introduced by Harvey and Lawson 
\cite{MR666108}
as a key example of a {\em calibrated geometry}: such submanifolds have
a local volume-minimizing property.  Unfortunately, very few examples
of special Lagrangian submanifolds are known in the compact case.

Let $X$ be a compact Calabi--Yau manifold.  
Physical arguments predict that for most such $X$, there should exist
a mirror partner $Y$, which is another compact Calabi--Yau manifold
whose physical 
(but not geometrical) properties are closely related to those of 
$X$.\footnote{For a brief review of
Calabi--Yau manifolds and mirror symmetry, see \cite{geomaspects}.  
More extensive
reviews can be found in \cite{MR1677117,Gross:LN,MR2003030}.}
However, as has been recognized since the
early days of mirror symmetry 
\cite{Ceresole:1992su,2param1,Hosono:1993qy,catp}, 
$X$ is expected to have a mirror partner
only when the complex structure of $X$ is sufficiently close to
a ``large complex structure limit point,'' which is a class of boundary
points in the compactified moduli space $\mathcal M_X$ characterized
by having maximally unipotent monodromy 
\cite{mirrorguide,delignelimit,compact}.  Strominger, Yau, and Zaslow
\cite{syz} argued on physical grounds that any such  
Calabi--Yau manifold should have a map
$\pi: X\to B$ whose general fiber $\pi^{-1}(b)$ is a special
Langrangian $n$-torus $T^n$; such a structure is called a 
{\em supersymmetric torus fibration of $X$}, since the special Lagrangian
condition is the geometric counterpart to the preservation of
half of the supersymmetry in a physical model.  Strominger, Yau, and
Zaslow also proposed that the mirror partner should be given,
to first approximation, by a compactification of the family of dual tori 
$\bigcup (\pi^{-1}(b))^\vee$.

It is worth saying a few words about the physical construction of
the mirror partner.  Strominger, Yau, and Zaslow argue that the
original Calabi--Yau manifold and its mirror partner should be related
by a physical construction known as ``T-duality on the $n$-torus
fibers.''  In the absence of holomorphic disks with boundaries on
the $n$-tori, this T-duality simply replaces each nonsingular torus by its 
dual torus
(cf.~\cite{MR2154821}), while doing something unknown at the 
singular fibers.  This description is expected to be modified when
holomorphic disks are present, but the precise effect of the holomorphic disks
 has not yet been worked out.  And as we shall
see, the current expectation (at least when $n=3$)
is that such holomorphic disks will be
present for at least some of the tori in the torus fibration.
This has made it difficult to formulate a mathematical version of
the original
Strominger--Yau--Zaslow proposal which is both precise and accurate.

Because the arguments used by Strominger, Yau, and Zaslow implicitly
assume that the Calabi--Yau metric is uniformly large, we put that
hypothesis in the following version of their existence conjecture.

\begin{conjecture}[Existence; cf.~\cite{syz}]
For any Calabi--Yau metric on a compact complex manifold $X$ of
complex dimension $n$ whose complex
structure is sufficiently close to a large complex structure limit point
and whose K\"ahler class is sufficiently 
deep in the K\"ahler cone,
there exists a supersymmetric torus fibration $\pi: X\to B$, where $B$
is a homology $n$-sphere.
\end{conjecture}

There is by now considerable indirect evidence in favor of this
conjecture, including an explicit construction in a (slightly degenerate)
limiting case \cite{GrossWilson:3tori}, as well as two strategies
(\cite{Ruan:quinticI,Ruan:quinticII,Ruan:quinticIII,Ruan:different} and 
\cite{Ruan:hypersurfacesI,Ruan:hypersurfacesII,Ruan:hypersurfacesIII})
for constructing weak forms of these fibrations for a certain
class of compact Calabi--Yau threefolds.
However, it has become clear that
proving this conjecture will require developing new techniques for
studying special Lagrangian submanifolds of a Calabi--Yau manifold.
In spite of our lack of tools to prove the conjecture, though, 
many qualitative features of supersymmetric torus
fibrations
have been inferred in various ways, and this paper is devoted to
explaining our best 
current (conjectural) understanding of those qualitative features.

We introduce the following terminology and notation.
Given a supersymmetric torus fibration $\pi: X\to B$, we let 
$\Sigma\subset X$ be the set of {\em singular points of fibers}\/ of $\pi$,
and let $\Delta=\pi(\Sigma)$ be the {\em discriminant locus}\/ 
of the fibration.

\section{Examples in Low Dimension}

In low dimension, 
supersymmetric torus fibrations of compact
Calabi--Yau manifolds are completely understood.

In the case of elliptic curves ($n=1$), any Calabi--Yau
 metric is flat and a special
Lagrangian $1$-torus is just a closed geodesic.  As is well known, if the
homology class is fixed then there is a fibration of the elliptic curve
over $B=S^1$ by closed geodesics in the specified class, with no singular
fibers.

In the case of K3 surfaces ($n=2$), an analysis is possible due to
the non-uniqueness of the compatible complex structure.  In fact, the
original paper of Harvey and Lawson \cite{MR666108}
showed that if $L\subset X$ is special Lagrangian, then there is
a different complex structure on $X$ (compatible with the given 
Calabi--Yau metric) such that $L\subset X$ is a complex submanifold.
Thus, a special Lagrangian
$T^2$ fibration can be interpreted in another complex structure
 as a holomorphic elliptic
fibration  $\pi:X\to B$ (with base $B=\mathbb{CP}^1=S^2$), 
and the structure of these is known in detail.
(In fact, a generic Ricci-flat metric
on a K3 surface admits such a fibration \cite{underlying}, so the existence
conjecture holds in this case.)  
Thanks to work of Kodaira \cite{MR0184257}, a complete classification of
possible singular fibers of such fibrations
is known: they are characterized by the conjugacy
class of the monodromy action on $H^1(T^2,\mathbb{Z})$.  The simplest
fibers, called {\it semi\-stable}, are associated to unipotent monodromy
transformations.  In an appropriate basis, the monodromy matrix takes the
form
\begin{equation} \label{eq:mono}
M=\left(\begin{array}{cc} 1&k\\0&1 \end{array} \right) \thinspace .
\end{equation}

The topology of a semi\-stable degeneration with $k=1$ is very familiar.
One of the cycles on the two-torus extends over the degeneration, and the
other ``vanishing'' cycle shrinks to a point; 
if we follow the torus around a loop encircling
the degeneration point in the base, there is a Dehn twist along the
vanishing cycle.  In spite of this twisting of the topology of the torus,
though, the total space of the fibration is non-singular.
(In the physics literature, the corresponding geometry is known as the
``Taub-NUT metric.'')

For the generic elliptic fibration of a K3 surface, all fibers are
semi\-stable, and there are exactly $24$
of them, each with $k=1$.  The monodromy data for such a generic
fibration (choosing an arbitrary base point $b\in B$)
gives a natural homomorphism
\begin{equation}\label{eq:22}
 \pi_1(B-\{P_1,\dots,P_{24}\},b)\to SL(2,\mathbb Z)
\end{equation}
whose generating loops all map to matrices conjugate to eq.~\eqref{eq:mono}.

The mirror partner of a given K3 surface (with a fixed Ricci-flat metric)
is known to be
another K3 surface with a different Ricci-flat metric \cite{stringK3}.  
When passing to the mirror, the monodromy
matrices $M$ are replaced by ${}^t\!M^{-1}$; since ${}^t\!M^{-1}$ 
is conjugate to
$M$, the monodromy data does not change.  
We will conjecturally extend this kind of ``topological''  mirror symmetry
statement
to dimension $3$ in the next section.
Note that, in any dimension, if we replace all nonsingular tori by their
dual tori, the monodromy matrices change as
$M\mapsto {}^t\!M^{-1}$.

\section{Smooth $T^3$ fibrations} \label{sec:smooth}

The first step in studying supersymmetric $T^3$ fibrations of
compact 
Calabi--Yau threefolds is to study more general $T^3$ fibrations, without
imposing the ``special Lagrangian'' condition.  
For Calabi--Yau hypersurfaces in toric varieties (of arbitrary
dimension), Zharkov \cite{Zharkov:torus} constructed a (topological)
$T^n$ fibration.
A general program to
understand such fibrations $\pi:X\to B$ (in dimension $3$)
for which the map $\pi$ is smooth
(i.e., $C^\infty$) was initiated by Gross 
\cite{Gross:slagI,Gross:slagII,Gross:Tmir}, and
parallel results were obtained by Ruan 
\cite{Ruan:quinticI,Ruan:quinticII,Ruan:quinticIII} in his study of 
a specific class of $T^3$ fibrations of Calabi--Yau hypersurfaces.
The monodromy of such fibrations and the topology of the singular
fibers was determined under a suitable assumption of genericity,
analogous to the assumption of ``generic elliptic fibration'' in the
case of K3 surfaces which guaranteed that all fibers were semi\-stable
with $k=1$.  We can summarize the analysis in a general conjecture
(which conjecturally extends their results to the general case).

\begin{conjecture}[Topology; cf.~\cite{Gross:slagI,Gross:slagII,Gross:Tmir,%
Ruan:quinticI,Ruan:quinticII,Ruan:quinticIII}] \label{con:topology}
Let $\pi: X\to B$ be a 
smooth
$T^3$
fibration\footnote{We stress that it is the map $\pi$ which is
smooth, not the fibers of the fibration.} of a compact 
Calabi--Yau threefold 
which is generic in a suitable sense.
Then
\begin{itemize}

\item[a)] The discriminant locus of the fibration
 is a trivalent graph $\Gamma$.

\item[b)] The topology near the edges of $\Gamma$ is modeled by the product 
of a cylinder with
the $k=1$
sem\-istable degeneration of two-tori.
In particular, a Dehn twist along the vanishing
cycle and a nonsingular total space are features of this topology.

\item[c)] For any loop around an edge of $\Gamma$, the monodromy on 
either $H^1\cong H_2$ or $H_1$ of the $3$-tori
is conjugate to
\begin{equation} \label{eq:mono3}
M=\left(\begin{array}{ccc} 1&0&1\\0&1&0\\0&0&1 \end{array} \right) \thinspace.
\end{equation}
In particular, both monodromy actions have a $2$-dimensional fixed plane.

\item[d)] The vertices of $\Gamma$ come in two types: near a 
{\em positive vertex},
the three monodromy actions on $H^1\cong H_2$ near the vertex have 
fixed planes whose intersection is $1$-dimensional, 
while the three monodromy actions
on $H_1$ have a common $2$-dimensional fixed plane.
In an appropriate basis, the monodromy matrices on $H_1$ take the form
\begin{equation} \label{eq:posmat}
\left(\begin{array}{ccc} 1&0&1\\0&1&0\\0&0&1 \end{array} \right) ,
\quad
\left(\begin{array}{ccc} 1&0&0\\0&1&1\\0&0&1 \end{array} \right) ,
\quad
\left(\begin{array}{ccc} 1&0&-1\\0&1&-1\\0&0&\hphantom{-}1 \end{array} 
\right) .
\end{equation}
On the other hand, near a {\em negative vertex},
the three monodromy actions on $H^1\cong H_2$ near the vertex have a 
common $2$-dimensional fixed plane, while the three monodromy actions
on $H_1$ have fixed planes whose intersection is $1$-dimensional.
In an appropriate basis, the monodromy matrices on $H_1$ take the form
\begin{equation} \label{eq:negmat}
\left(\begin{array}{ccc} 1&0&1\\0&1&0\\0&0&1 \end{array} \right) ,
\quad
\left(\begin{array}{ccc} 1&1&0\\0&1&0\\0&0&1 \end{array} \right) ,
\quad
\left(\begin{array}{ccc} 1&-1&-1\\0&\hphantom{-}1&\hphantom{-}0\\0
&\hphantom{-}0&\hphantom{-}1 \end{array} \right) .
\end{equation}

\item[e)] 
The fiber of 
$\pi$
over any point of $B$ other than a vertex of $\Gamma$ has a fixed point free
$U(1)$ action and in particular has Euler characteristic $0$.  The
fiber of 
$\pi$
 over each positive vertex  has Euler characteristic $1$,
and the fiber over each negative vertex
has Euler characteristic $-1$.
(In fact, Gross \cite{Gross:Tmir} and Ruan \cite{Ruan:quinticIII}
gave explicit descriptions of these singular fibers, but we will
not reproduce those descriptions here.)
\end{itemize}
\end{conjecture}

There is an induced global monodromy action on $H_1$
\begin{equation}\label{eq:32}
 \pi_1(B-\Delta,b) \to SL(3,\mathbb Z),
\end{equation}
whose generators satisfy the conditions spelled out in the conjecture.
If we have a compact Calabi--Yau threefold and its mirror partner, with smooth
$T^3$ fibrations whose nonsingular fibers are dual to each other, then the
monodromy transformations will be related as
$M\mapsto {}^tM^{-1}$.  This implies
 that the r\^oles of the positive and negative vertices in the 
fibration  are reversed between a Calabi--Yau threefold and its mirror
partner.
Since part (e) of the conjecture implies that the topological
Euler number of $X$ can be calculated via
\begin{equation}\label{eq:33}
 \chi_{top}(X)=\#\{\text{positive vertices}\} - \#\{\text{negative vertices}\},
\end{equation}
the effect of mirror symmetry on monodromy then shows that the Euler
number changes sign:
\begin{equation}\label{eq:34}
 \chi_{top}(Y)
=\#\{\text{negative vertices}\} - \#\{\text{positive vertices}\}
=-\chi_{top}(X),
\end{equation}
as expected from physical mirror symmetry arguments.

\begin{figure}
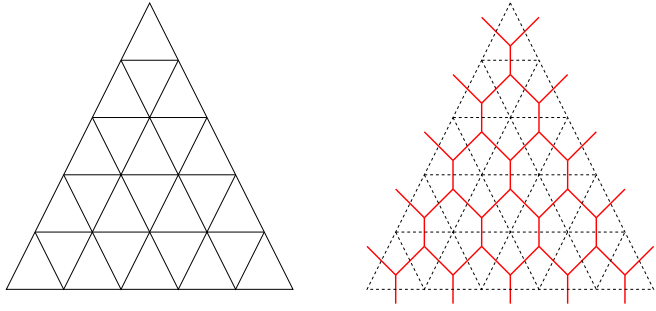

\begin{center}
\includegraphics[scale=0.3]{5trianglea.mps}
\qquad
\includegraphics[scale=0.3]{5triangle3b.mps}
\end{center}
\caption{A triangulation of a two-dimensional face of
the Newton polytope of the quintic threefold, and
the corresponding dual graph.}\label{fig:triangulation}
\end{figure}

In \cite{Gross:Tmir}, Gross showed how to go
further, and use the data from a generic smooth $T^3$ fibration of
a given compact Calabi--Yau threefold to construct a manifold which is a 
candidate mirror partner.  The transpose inverse of the original
monodromy representation produces a mirror monodromy representation,
describing the monodromy on the family of dual tori (with positive
and negative vertices reversed).
Gross proved a ``Reconstruction Theorem:''
the  family of dual tori can be
completed to a compact topological
manifold with a smooth $T^3$ fibration, satisfying the
properties stated in the conjecture.

\section{Combinatorics of $\Gamma$}

\begin{figure}
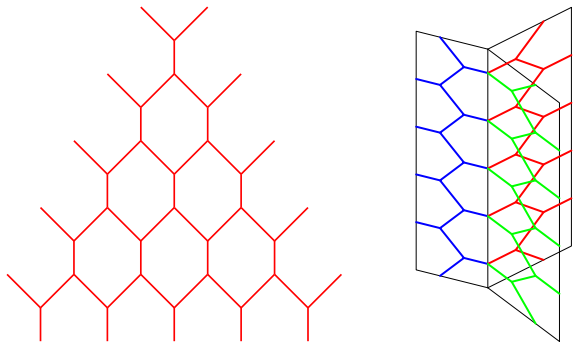

\begin{center}
\includegraphics[scale=0.35]{5triangle4a.mps}
\qquad
\includegraphics[scale=0.25]{junction.mps}
\qquad
\hphantom{.}
\end{center}
\caption{Portions of $\Gamma$ for the quintic threefold.}
\label{fig:graphpiece}
\end{figure}

For smooth $T^3$ fibrations of a compact Calabi--Yau threefold,
the combinatorics of the graph $\Gamma$ are beautiful and intricate.
For example, in the case of a quintic hypersurface in $\mathbb {CP}^3$,
the graph depends on  choices of triangulations of the two-dimensional
faces of the Newton polytope of the defining equation; one such 
choice is shown on the left side of Figure~\ref{fig:triangulation}.
One constructs the dual graph 
of each such triangulation, in which each face of the triangulation gives
a vertex of the graph, and each edge of the triangulation is crossed
by an edge of the graph (as shown on the right side of
Figure~\ref{fig:triangulation}).  That
dual graph then becomes a piece of $\Gamma$ (illustrated on
the left side of Figure~\ref{fig:graphpiece}) in which each vertex
is ``negative.''
The pieces are assembled according
to the combinatorics of the Newton polytope, in which the free
ends of the dual graph meet free ends from other faces of the
Newton polytope, forming trivalent vertices which are the
``positive'' vertices of $\Gamma$.  The pieces thus attach three at a time;
a neighborhood of one such attachment is illustrated on
the right side of Figure~\ref{fig:graphpiece}

A general description of these graphs, for Calabi--Yau complete intersections
in toric varieties, was given by 
Haase and Zharkov \cite{HZ:T3I,HZ:T3II,HZ:T3III} and by Gross \cite{Gross:BB}.
Their results\footnote{The proofs of some of the results stated by
Gross \cite{Gross:BB} were deferred to another paper which has not yet
appeared.} show that the classes of mirror pairs described by
Batyrev \cite{MR1269718}
and by Batyrev and Borisov \cite{borisov-ms,MR1463173} admit smooth
$T^3$ fibrations satisfying Conjecture~\ref{con:topology}, and that those
$T^3$ fibrations are mirror duals of each other.

\begin{figure}
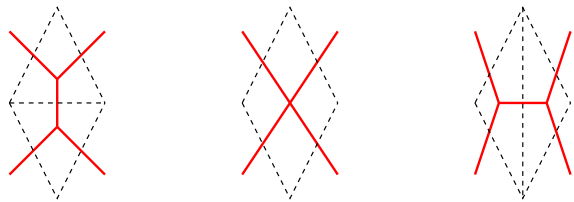

\begin{center}
\includegraphics[scale=0.5]{flop+.mps}
\qquad\qquad
\includegraphics[scale=0.5]{flopsing.mps}
\qquad\qquad
\includegraphics[scale=0.5]{flop-.mps}
\end{center}
\caption{The change of triangulation corresponding to a flop.}
\label{fig:flop}
\end{figure}

This theory of smooth $T^3$ fibrations of compact Calabi--Yau threefolds
can also be related to some  topological aspects of mirror
symmetry which have played important r\^oles 
in the physics literature \cite{Witten:1993yc,catp,bhole}.
First, the construction of $\Gamma$ for a Calabi--Yau hypersurface $X$
in a toric
fourfold depends on a choice of triangulation
of the faces of the Newton polytope, and  this choice of triangulation is
equivalent to a choice of large complex structure limit point in the moduli
space \cite{GKZ2,catp}.  Mirror symmetry (as developed in the physics
literature) offers an alternate interpretation:
each large complex structure limit point corresponds to a different
birational model of the mirror partner, and the K\"ahler cones of the
birational models fit together into a common space (after complexification),
mirroring the complex structure moduli space $\mathcal{M}_X$ \cite{beyond}.
In this interpretation,
the choice of birational model depends explicitly on a choice of
triangulation of the Newton polytope (which describes the mirror toric
fourfold in Batyrev's construction \cite{MR1269718}).  The simplest
birational change---a ``flop''---is realized by the simple change of
triangulation
illustrated in Figure~\ref{fig:flop}.  As Gross pointed out
\cite[Remark 4.5]{Gross:Tmir}, the corresponding change in dual graph
(also illustrated in the figure) has the correct monodromy properties
to be allowed as a new graph $\Gamma'$.
One expects that by appropriately 
varying the complex structure and the K\"ahler metric,
the connecting edge in the original graph will shrink to zero length, giving
the $4$-valent vertex shown in the intermediate stage; further variation
then causes a new connecting edge to grow, changing the topology to that
of $\Gamma'$.

\begin{figure}
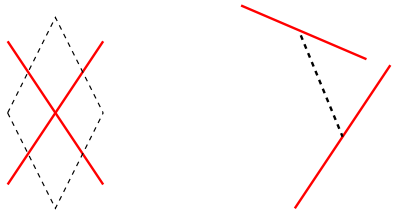

\begin{center}
\includegraphics[scale=0.5]{flopsing.mps}
\qquad\qquad
\includegraphics[scale=0.5]{conifold.mps}
\end{center}
\caption{Expected change of graph for a conifold transition.}
\label{fig:conifold}
\end{figure}

The intermediate step illustrated in the middle of Figure~\ref{fig:flop}
is a partial triangulation, which is expected to correspond to a 
``conifold'' singularity on the mirror Calabi--Yau threefold.  A second
topological feature of mirror symmetry is the ``conifold transition''
\cite{Candelas:1987kf,bhole}
in which that conifold singularity is resolved with a small blowup,
rather than being smoothed with a change of complex structure.
(The conifold singularity is on the mirror partner, but this transition
can also be described on the original Calabi--Yau threefold
\cite{looking}.)
The graph $\Gamma$ appears to change as follows, as proposed independently
by Gross \cite{Gross:slagEx} and Ruan \cite{Ruan:lagmir}: after shrinking
the connecting edge to zero size, leaving a $4$-valent vertex, the two
arms of the graph crossing at that vertex are separated into different
planes, as illustrated in Figure~\ref{fig:conifold}.  As Gross and
Ruan verify, this change is compatible with the monodromies around the edges
and produces the expected change in topological Euler characteristic
for a conifold transition.  However, the relation between this construction
and more global versions of the conifold transition remains mysterious
and needs further study.

\section{Affine structures on the base} \label{sec:affine}

Local moduli for a compact
special Lagrangian submanifold $L$ of a compact Calabi--Yau manifold
$X$ were determined by McLean \cite{MR1664890}: the deformation space is
smooth, and its tangent space is canonically identified with the
space of harmonic $1$-forms of $L$.
Hitchin \cite{dg-ga/9711002} used this identification
to construct two affine structures on $B-\Delta$
(if $\pi$ is smooth), which geometrize the monodromy transformations
that occurred in Conjecture~\ref{con:topology}.

If $V$ is a normal vector field to $L$ in $X$, then the contraction
$\iota(V)\omega$ of $V$ with the K\"ahler form gives a harmonic
$1$-form on $L$, and the contraction $\iota(V)\Omega$ of $V$ with
the holomorphic $n$-form gives a harmonic $(n-1)$-form on $L$.  These
constructions give rise to the affine structures, which can be seen 
by considering periods of the harmonic
 forms.  If $\frac{\partial}{\partial t_i}$, 
\dots, $\frac{\partial}{\partial t_n}$ are vector fields on the
deformation space $\mathcal{M}$ which span the tangent space to 
$\mathcal{M}$ at $[L]$,
and $A_1$, \dots, $A_n$ is a basis for
$H_1(L,\mathbb Z)$, then we can form the period matrix
\begin{equation}
 \lambda_{ij}= \int_{A_i} \iota(\frac{\partial}{\partial t_j})\omega.
\end{equation}
The $1$-forms $\sum_j \lambda_{ij}dt_j$ on the deformation space 
$\mathcal{M}$
are closed, and can be integrated to give local coordinates
$u_1$, \dots, $u_n$ on $\mathcal{M}$ at $[L]$
satisfying $du_i=\sum_j \lambda_{ij}dt_j$.  
Such coordinate
systems provide an affine structure, that is, the transition functions
between any two such coordinate systems lie in the affine group
$\mathbb{R}^n\rtimes\operatorname{GL}(\mathbb{Z}^n)$.
Intrinsically, the lattice $\mathbb{Z}^n$ should be identified with
$H_1(L,\mathbb Z)$ in this case, and this affine structure carries the
information about the monodromy on $H_1(L,\mathbb Z)$.

Similarly, if $B_1$, \dots, $B_n$ is a basis for $H_{n-1}(L,\mathbb Z)$,
then we can form the period matrix
\begin{equation}
 \mu_{ij}= \int_{B_i} \iota(\frac{\partial}{\partial t_j})\Omega.
\end{equation}
The $1$-forms $\sum_j \mu_{ij}dt_j$ on the deformation space 
$\mathcal{M}$
of $L$ are closed, and can be integrated to give local coordinates
$v_1$, \dots, $v_n$ on 
$\mathcal{M}$ satisfying $dv_i=\sum_j \mu_{ij}dt_j$.  
Such coordinate
systems provide the other affine structure, which carries the information about
the monodromy on $H_{n-1}(L,\mathbb Z)\cong H^1(L,\mathbb Z)$.

Hitchin shows that these two affine structures are related by a Legendre
transform with respect to a suitable locally defined function $K$
(which also determines a canonical metric on the deformation space).  This is a
version of mirror symmetry which is formulated
strictly on the base of the torus fibrations,
a notion which was further developed in \cite{GrossWilson:LCSL,KS:mirror}.

\section{The large complex structure limit}

When $X$ is a compact Calabi--Yau threefold, it is 
expected that the map $\pi$ giving
a supersymmetric torus fibration $\pi:X \to B$ will be only
 piecewise smooth, so the analysis of Sections~\ref{sec:smooth} 
and \ref{sec:affine} does
not directly apply.  
Gross showed \cite{Gross:slagII} that if $\pi$ is smooth, the
discriminant locus $\Delta\subset B$ must have codimension two; in
the more general ``piecewise smooth'' case, $\Delta$
may have codimension $1$. However, we do expect that
$\Delta$ will always have a retraction onto a subset $\Gamma$ of
codimension two.  And there is a particular
limiting situation in which
this retraction should become evident: the large complex structure
limit.

In fact, Gross--Wilson \cite{GrossWilson:LCSL}
and Kontsevich--Soibelman \cite{KS:mirror} have formulated a 
precise conjecture
about the large complex structure limit of a supersymmetric torus
fibration.\footnote{Note there has been substantial additional progress on 
these limits and on various
structures on the base in subsequent work of Kontsevich
and Soibelman \cite{KS:affine,arXiv:0811.2435}, reviewed elsewhere
in this volume.}

\begin{conjecture}[Large Complex Structure Limit; 
cf.~\cite{GrossWilson:LCSL,KS:mirror}]
Let $\mathcal{X} \to S$ be a maximally unipotent degeneration of
compact simply-connected Calabi--Yau manifolds of complex dimension $n$,
degenerating at $0\in S$,
let $s_i\in S$ be a sequence with $\lim s_i=0$, and let $g_i$ be
a sequence of Ricci-flat metrics on $X_{s_i}$ with 
diameter bounded above and below.  
Then there exists a subsequence $(X_{s_{i_j}},g_{i_j})$
which converges in the sense of
Gromov--Hausdorff \cite{MR682063} to a metric space
$(X_\infty, d_\infty)$, where $X_\infty$ is homeomorphic to
to the sphere
$S^n$, and $d_\infty$ is induced by a Riemannian metric on
$X_\infty - \Gamma_\infty$ for some $\Gamma_\infty\subset X_\infty$
of codimension two.
\end{conjecture}

Following the affine structure to this limit, one expects to find a
limiting affine structure, and in fact the discriminants $\Delta_{i_k}$
should have collapsed to $\Gamma_\infty$ in the limit.

\section{Non-compact examples of special Lagrangian fibrations} 
\label{sec:noncompact}

Harvey and Lawson's original paper about calibrations \cite{MR666108} gave
an explicit example of a special Lagrangian 
fibration.\footnote{Further development of examples of this type was made by
Goldstein \cite{MR1865243} and Gross \cite{Gross:slagEx}, 
and they have been extensively
studied in the physics literature (e.g., in
\cite{agva,akv,Aganagic:2003db}).}
Define $f:\mathbb{C}^3\to\mathbb{R}^3$ by
\begin{equation}
 f(z_1,z_2,z_3)
=(\operatorname{Im}(z_1z_2z_3),|z_1|^2-|z_2|^2,|z_1|^2-|z_3|^2).
\end{equation}
Then the fibers of $f$ are special Lagrangian, and are all invariant
under the action of the diagonal torus with determinant $1$:
\begin{equation}
\{\operatorname{diag}(e^{i\theta_1},e^{i\theta_2},e^{i\theta_3})\ |\
\theta_1+\theta_2+\theta_3=0\}.
\end{equation}

\begin{figure}
\begin{center}
\includegraphics[scale=0.5]{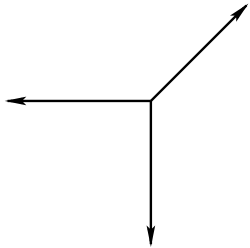}
\end{center}
\caption{The discriminant locus in the $(x_2,x_3)$-plane
for the Harvey--Lawson
fibration of $\mathbb{C}^3$.}
\label{fig:trivalent}
\end{figure}

The singularities of fibers are located where $z_i=z_j=0$ for some pair of
indices $i$ and $j$; the images of these give three rays within the
plane $\{x_1=0\}\subset \mathbb{R}^3$, namely
(i) $x_2=0$, $x_3\le0$, (ii) $x_2\le0$, $x_3=0$, and (iii) $x_2=x_3\ge0$.
This is illustrated in Figure~\ref{fig:trivalent}.

The nonsingular fibers are all homeomorphic to $T^2\times\mathbb{R}$.
Note that whenever $x_1=0$, if we write $z_j=r_je^{i\theta_j}$ then
either $\theta_1+\theta_2+\theta_3=0$ (in which case 
$\operatorname{Re}(z_1z_2z_3)\ge0$), or 
$\theta_1+\theta_2+\theta_3=\pi$ (in which case
$\operatorname{Re}(z_1z_2z_3)\le0$).  Thus, there are two natural
subsets $f^{-1}(0,x_2,x_3)^\pm$ of the fiber $f^{-1}(0,x_2,x_3)$,
distinguished by the sign of $\operatorname{Re}(z_1z_2z_3)$.
These subsets meet along 
\begin{equation}
f^{-1}(0,x_2,x_3)\cap\{z_1z_2z_3=0\}.
\end{equation}

When the fiber is smooth, each subset is a manifold with boundary,
and they meet along their common boundary.  However, when the fiber
is singular and $(x_2,x_3)\ne(0,0)$, each subset is itself a
smooth special Lagrangian submanifold, homeomorphic to
$S^1\times\mathbb{R}^2$.  Note that $f^{-1}(0,0,0)$ is
a special case: as Harvey and Lawson pointed out, each subset
$f^{-1}(0,0,0)^\pm$ is a cone over $T^2$, and those cones meet
precisely at the origin in $\mathbb{C}^3$.

Joyce \cite{Joyce:SYZ}
builds some new special Lagrangian fibrations by carefully combining
subsets of Harvey--Lawson fibers.  Let
\begin{equation}
\begin{aligned} N_a^+ = \{&|z_1|^2-a = |z_2|^2+a
= |z_3|^2+|a|,\\ &\operatorname{Im}(z_1z_2z_3)=0,
\operatorname{Re}(z_1z_2z_3)\ge0 \}
\end{aligned}
\end{equation}
and
\begin{equation}
\begin{aligned} N_a^- = \{&|z_1|^2-a = |z_2|^2+a
= |z_3|^2+|a|,\\ &\operatorname{Im}(z_1z_2z_3)=0,
\operatorname{Re}(z_1z_2z_3)\le0 \}.
\end{aligned}
\end{equation}
Then 
\begin{equation}
N_a^+=\begin{cases} f^{-1}(0,2a,2a)^+ & \text{when } a\ge0,\\
f^{-1}(0,a,0)^+ & \text{when } a\le 0.
\end{cases}
\end{equation}
and similarly for $N_a^-$.

To build a special Lagrangian fibration, Joyce considers translations of these
manifolds for $c\in\mathbb{C}$.  Let
\begin{equation}
\begin{aligned} N_{a,c}^\pm = \{&|z_1|^2-a = |z_2|^2+a
= |z_3-c|^2+|a|,\\ &\operatorname{Im}(z_1z_2(z_3-c))=0,
\pm\operatorname{Re}(z_1z_2(z_3-c))\ge0. \}
\end{aligned}
\end{equation}
These can be made the fibers of special Lagrangian fibrations by
defining 
$F^\pm:\mathbb C^3\to \mathbb R\times\mathbb C$ by
\begin{equation}
F^\pm(z_1,z_2,z_3)=
\begin{cases} (\frac12(|z_1|^2-|z_2|^2), z_3) & \text{if } z_1=z_2=0 \\ 
(\frac12(|z_1|^2-|z_2|^2), z_3\mp\frac{\bar{z}_1\bar{z}_2}{|z_1|}) 
& \text{if } |z_2|^2 \le |z_1|^2 \ne 0 \\
(\frac12(|z_1|^2-|z_2|^2), z_3\mp\frac{\bar{z}_1\bar{z}_2}{|z_2|}) 
& \text{if } |z_2|^2 > |z_1|^2
 \end{cases}.
\end{equation}
With this definition, 
$(F^{\pm})^{-1}(a,c)=N^\pm_{a,c}$.

Notice that the fibrations $F^\pm$ are only piecewise smooth, 
and that the discriminant locus in each case is $\{(0,c)\}\subset
\mathbb{R}\times\mathbb{C}$, which has codimension $1$.  Notice
also that when $a>0$, both $N^+_{a,c}$ and $N^-_{a,c}$
contain the boundary of the
holomorphic disk $\{|z_1|^2\le a\}$ of area $2\pi a$, and that
when $a<0$, both $N^+_{a,c}$ and $N^-_{a,c}$
contain the boundary of the
holomorphic disk $\{|z_2|^2\le -a\}$ of area $-2\pi a$.
In some sense, these shrinking disks are ``responsible'' for the
singularity being created at $a=0$.

\begin{figure}
\begin{center}
\includegraphics[scale=1]{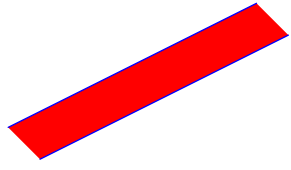}
\end{center}
\caption{The discriminant locus as a ribbon.}\label{fig:ribbon}
\end{figure}

To go further, Joyce invokes the extensive theory which he
developed in \cite{Joyce:U1I,Joyce:U1II,Joyce:U1III} concerning the structure
of special Lagrangian $3$-manifolds with a $U(1)$ action.  Using that
theory, he is able to construct \cite[Theorem 6.5]{Joyce:SYZ}
a special Lagrangian fibration 
$\widehat{F}:\mathbb{C}^3\to\mathbb{R}^3$ 
whose discriminant locus is a {\em ribbon}, that
is, the locus $\{(0,x_2,x_3)\ |\ 0\le x_2\le 1\}\subset\mathbb{R}^3$,
as illustrated in Figure~\ref{fig:ribbon}.
There are several important properties of this example of Joyce's.
First, the fiber over an interior point of the ribbon has two 
singularities---one locally modeled by $F^+$ and the other locally modeled
by $F^-$.  Second, as in those local models, there are holomorphic disks
with boundary in the fiber for $x_1\ne0$ (and $0<x_2<1$),
whose area approaches $0$ as $x_1$ approaches $0$; in fact, there
is one such holomorphic disk for each of the two singular points.

\begin{figure}
\begin{center}
\includegraphics[scale=0.8]{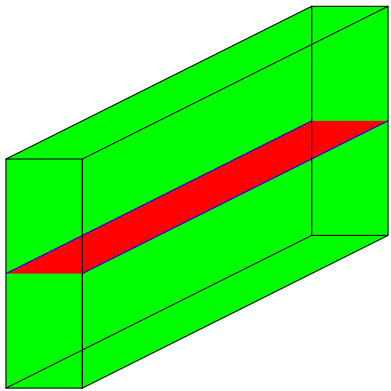}
\end{center}
\caption{The nearby tori which bound a holomorphic disk.}
\label{fig:thickribbon}
\end{figure}

Third, as we approach the boundary of the strip  within the plane $x_1=0$, 
something interesting
happens: the two bounding circles approach each other and the holomorphic
disks cancel out as the boundary of the strip (either $x_2=0$ or $x_2=1$)
is reached.  There are no holomorphic disks when $x_2<0$ or $x_2>1$.
The region in which holomorphic disks are present is illustrated in
Figure~\ref{fig:thickribbon}.

Thus, along the boundary of the strip, one has a singularity
of multiplicity $2$ (in an appropriate sense), which bifurcates into a
pair of singularities in the middle of the strip, and those singularities
rejoin at the other boundary.

Note that the plane which contains the discriminant locus can be
identified intrinsically using local affine coordinates.  The cycle
$\gamma\in H_1(N_{a,c}^\pm,\mathbb{Z})$ which bounds a holomorphic disk is the
vanishing cycle for the family, and it defines a dual subspace 
in $\gamma^\perp\subset H^1_c(N_{a,c}^\pm,\mathbb{Z})$, which can be locally
identified with the plane containing the discriminant locus.
This plane can also be characterized as 
the monodromy-invariant plane in the compactly-supported
cohomology of the fiber.

Joyce conjectures that his examples exhibit generic behavior.
In fact, even Lagrangian fibrations (not just special
Lagrangian fibrations) are expected to
exhibit these phenomena, as explained in \cite{math.SG/0611139}.

\begin{conjecture}[Singular Fibers; cf. \cite{Joyce:SYZ}]
\label{con:joyce}
Let $\pi:X\to B$ be a supersymmetric $T^3$ fibration of a compact Calabi--Yau
threefold with respect to a Calabi--Yau metric
whose compatible complex structure is sufficiently close to a
large complex structure limit point, and whose K\"ahler class is
sufficiently deep in the K\"ahler cone.  Then 
\begin{enumerate}
\item The discriminant locus $\Delta\subset B$ has codimension one.
In affine coordinates, $\Delta$
 is locally contained in the plane corresponding to 
the monodromy-invariant subspace of $H^1(\pi^{-1}(b),\mathbb{Z})$
for $b$ near $\Delta$.
\item The fiber over the general point of $\Delta$ has two singular
points, one of which is modeled locally by $F^+$ and the other
of which is modeled locally by $F^-$.
\item The fiber over the general point of the boundary of $\Delta$
is modeled locally by $\widehat{F}$.
\item Let $\mathcal{H}\subset B$ be the set of fibers which contain
the boundary of at
least one holomorphic disk in $X$, or are singular.  Then the
boundary of $\Delta$ is contained in the boundary of $\mathcal H$.
\end{enumerate}
\end{conjecture}

The last statement about which fibers contain the boundaries of
 disks was not conjectured
by Joyce, but is consistent with the behavior exhibited by his example
$\widehat{F}$ (as illustrated in Figure~\ref{fig:thickribbon}).

\section{Amoebas}
\label{sec:amoebas}

A common feature of the constructions of Zharkov \cite{Zharkov:torus}
and Ruan \cite{Ruan:quinticI,Ruan:quinticII,Ruan:quinticIII,Ruan:different},
which dovetails nicely with the analysis of Joyce
described in Section~\ref{sec:noncompact}, is
the description of the discriminant locus of a supersymmetric torus
fibration as an {\em amoeba}.  Amoebas were introduced by
Gelfand, Kapranov, and Zelevinsky \cite{GKZ2}; 
we will briefly review the theory,
following 
Mikhalkin \cite{math.AG/0403015}
(see also 
\cite{Gross:amoebas}).

\begin{figure}
\begin{center}
\includegraphics[scale=0.5]{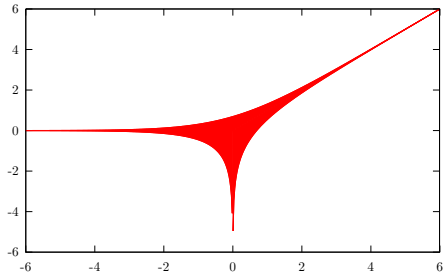}
\end{center}
\caption{
Amoeba of $z+w+1$.}
\label{fig:lineamoeba}
\end{figure}

\begin{figure}
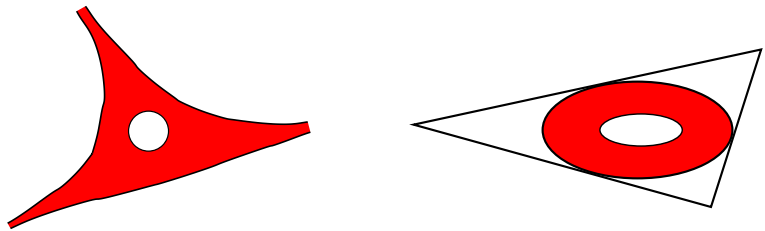

\begin{center}
\includegraphics[scale=.17]{am1a.mps}
\qquad
\includegraphics[scale=.3]{am2b.mps}
\end{center}
\caption{An amoeba with one hole, and its compactification.} \label{fig:amoeba}
\end{figure}

Let $f=\sum a_Ix^I$ be a Laurent polynomial in $n$ complex variables.
(Here, $I$ is a multi-index with negative powers allowed, but $f$ has only
finitely many non-zero terms.)  The {\em amoeba of $f$}\/ is the
set $\mathcal{A}_f=\operatorname{Log}(V_f)$, where
$V_f=\{\vec{z}\in (\mathbb{C}^*)^n\ |\ f(\vec{z})=0\}$,
and $\operatorname{Log}:(\mathbb{C}^*)^n\to \mathbb{R}^n$ is defined by
\begin{equation}
 \operatorname{Log}(z_1,\dots,z_n)=(\log |z_1|, \dots, \log |z_n|).
\end{equation}
A simple example is given by $f(z,w)=z+w+1$, which is often chosen because
it can be graphed exactly, as in Figure~\ref{fig:lineamoeba}.
More complicated examples have ``holes'' in the amoeba, as indicated
on the left side of Figure~\ref{fig:amoeba}.

The Laurent polynomial $f$ has an associated Newton polytope $\Delta_f$
and toric variety $T_f$, and there is a moment map
$\mu:T_f\to \Delta_f\subset\mathbb{R}^n$ once a symplectic structure
has been chosen on $T_f$.  The closure of $\mu(V_f)$
is called the {\em compactified amoeba of $f$}.  Note that 
$\operatorname{Log}$ and
$\mu$ are closely related: the interior of the image $\Delta_f$ of
$\mu$ is mapped homeomorphically by $\operatorname{Log}\circ \mu^{-1}$ to
all of $\mathbb{R}^n$.  An example of a compactified amoeba is
illustrated on the right side of Figure~\ref{fig:amoeba}.

Forsberg, Passare, and Tsikh \cite{MR1752241} showed that each component
of the complement
$\mathbb{R}^n-\mathcal{A}_f$ is convex, and that there is
an injective map from the set of components to the lattice points
$\Delta_f\cap\mathbb{Z}^n$ in the Newton polytope, in which the bounded
components (the ``holes'') map to points in the interior of $\Delta_f$.
For a given polyhedron $\Delta$, there exist functions 
$f$ with $\Delta_f=\Delta$ whose amoebas
have
the maximum number of holes, but typically there also exist functions 
with $\Delta_f=\Delta$
whose amoebas have
fewer holes.

In the case $n=2$,  a formula of Baker \cite{baker-newton,MR487230}
identifies the genus of a smooth compactification of the affine curve
$V_f$ with the number of interior lattice points in the Newton polytope.
Thus, in that case, the maximum number of holes coincides with the genus.
There is an associated topological picture when the number of
holes is maximal: the map $\operatorname{Log}$
will in this case be $2$-to-$1$ over the interior of the amoeba, and
$1$-to-$1$ on the boundary of the amoeba.  
It is easy to see that this gives the right
answer for the genus.

Additional information about the amoeba can be obtained by considering 
the {\em Ronkin function}\/ 
$N_f: \mathbb{R}^n \to \mathbb{R}$ defined by
\begin{equation}
 N_f(x_1,\dots,x_n) = \frac1{(2\pi i)^n}
\int_{\mu^{-1}(x_1,\dots,x_n)} \log | f(z_1,\dots,z_n)| \frac{dz_1}{z_1}
\wedge \dots \wedge \frac{dz_n}{z_n}.
\end{equation}
Ronkin \cite{MR1833516} and
Passare--Rullg{\aa}rd \cite{MR2040284} show that $N_f$ is well-defined
on all of $\mathbb{R}^n$, is convex over $\mathcal{A}_f$, and is locally
linear on the complement of $\mathcal{A}_f$.  Let $\{E\}$ be the set
of components of the complement, and let $N_E$ be the extension of
$N_f|_E$ to a linear function on all of $\mathbb{R}^n$.
Passare and Rullg{\aa}rd define
\begin{equation}
 N_f^\infty = \max_E N_E,
\end{equation}
which is a piecewise linear function on $\mathbb{R}^n$, and then define the
{\em spine of the amoeba $\mathcal{A}_f$} to be
the set $S_f\subset\mathbb{R}^n$ of points
at which  the function $N_f^\infty$ is not locally linear.
A key theorem of \cite{MR2040284} is that the spine $S_f$ is a strong
deformation retract of the amoeba $\mathcal{A}_f$.

Note that the spine of the amoeba which is shown in Figure~\ref{fig:lineamoeba}
is precisely given by Figure~\ref{fig:trivalent}.

Ruan \cite{Ruan:Newton} observed that for a Calabi--Yau hypersurface $X$ in a 
toric fourfold $T$ which is close to the large complex structure limit,
the intersections $C_{jk}=X\cap T_j\cap T_k$ with pairs of toric divisors
have amoebas with the maximum number of holes, and these amoebas retract
to their spines as the large complex structure limit is approached.
The spines in fact form the pieces of the graph $\Gamma$ used to describe
a topological $T^3$ fibration, which are dual graphs of appropriate 
triangulations of the Newton polytopes (as was illustrated in
Figure~\ref{fig:triangulation}).

\begin{figure}
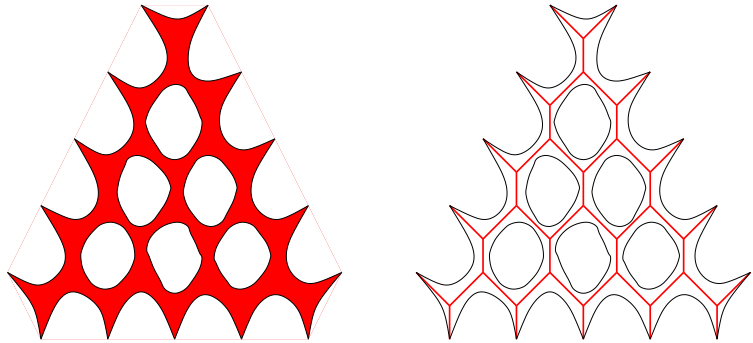

\begin{center}
\includegraphics[scale=0.35]{5triangle7.mps}
\qquad
\includegraphics[scale=0.35]{5triangle5a.mps}
\end{center}
\caption{Amoeba and spine for $C_{jk}$ of the quintic threefold.}
\label{fig:spine}
\end{figure}

An amoeba for $C_{jk}$, together with its spine, is shown in the case
of the quintic threefold in Figure~\ref{fig:spine}.  
The spine is precisely the graph
which occurred on the left side of Figure~\ref{fig:graphpiece}.

As in the case of $\Gamma$ itself,
moving among different large complex structure limit
 points causes the combinatorics of the triangulation to change 
(as discussed at the end of Section~\ref{sec:smooth});
we expect a corresponding change in the combinatorics
of the amoebas.

\section{Reconstruction in the Lagrangian case}

Gross's reconstruction theorem produced a topological $6$-manifold out
of the data describing a smooth $T^3$ fibration, and this provides
(in principle) a method for constructing mirror partners when they are
not known, provided that one has a smooth $T^3$ fibration.

Casta\~no-Bernard and Matessi 
\cite{math.SG/0611139} have proved an analogous theorem which produces
a compact symplectic $6$-manifold with a piecewise smooth Lagrangian $T^3$
fibration, starting from the data describing the smooth fibers of this
fibration.\footnote{Note that this is still not a {\em special}\/ Lagrangian
fibration, but provides an important intermediate step between the
cases of topological fibration and special Lagrangian fibration.}
The starting data in this case is the affine structure on the base---a
refinement of the simple monodromy data which Gross's theorem needed.
The key technique of stitching together smooth $T^3$ fibrations along
a common boundary had been developed in  earlier work of these authors
\cite{math.SG/0603498,MR2201567}.

A particularly interesting feature of Casta\~no-Bernard and Matessi's
construction is the behavior of the discriminant locus $\Delta$.
Their fibrations have a discriminant locus which is a trivalent graph
near the positive vertices, but has codimension $1$ near the negative
vertices.  The discriminant locus retracts onto a trivalent graph;
the inverse ``thickening'' of parts of this
 graph to a codimension $1$ set replaces
each neighborhood of a negative vertex with an amoeba-like shape which
retracts back to the graph.  Moreover, the fibration is smooth
outside of a set which retracts to (a subset of) the graph.
This result, when combined with the
discussion in Section~\ref{sec:amoebas}, helps to motivate our conjectures
in the next two Sections.

\section{The geometry of $T^3$ fibrations}

Prior to
Joyce's analysis of the structure of special Lagrangian fibrations
\cite{Joyce:SYZ},
there had been speculation that supersymmetric $T^3$ fibrations of compact
Calabi--Yau threefolds would always be smooth, so that the detailed
structure (in the generic case) would be the one given in 
Conjecture~\ref{con:topology}.  However, Joyce's analysis prompted many
of us to rethink the question, and to try to formulate properties
analogous to those of Conjecture~\ref{con:topology} which we would expect
supersymmetric $T^3$ fibrations to have.
Once such formulation appears in Conjecture~\ref{con:geometry} below.

\begin{figure}
\begin{center}
\includegraphics[scale=0.8]{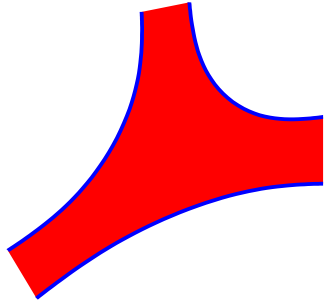}
\end{center}
\caption{The discriminant locus near a negative vertex.}\label{fig:negvertex}
\end{figure}

In Conjecture~\ref{con:topology}, the discriminant locus is a graph
$\Gamma$, but in general we should expect a discriminant locus $\Delta$
which only retracts to a graph $\Gamma$.
As Joyce pointed out 
(and was already mentioned in Conjecture~\ref{con:joyce}),
the first thing to expect is that the edges of the graph $\Gamma$
should thicken to ribbons; moreover, one
should see two singular points in each fiber over an interior
point of the ribbon, with the two points coming together to a single
singular point along the edges of the ribbon.
The next thing to expect
was also proposed by Joyce \cite{Joyce:SYZ}: 
since at a ``negative'' vertex, the local monodromy transformations 
share a common $2$-dimensional fixed plane, $\Delta$ should remain
planar, and the 
negative vertex should be replaced
by a ``trivalent ribbon'' 
of the sort illustrated in Figure~\ref{fig:negvertex}.
(This is the same structure found by Casta\~no-Bernard and Matessi
\cite{math.SG/0611139} in the Lagrangian case.)
Near a ``positive'' vertex, the three planes containing parts of $\Delta$
share a common line but are distinct; 
Joyce also made a specific proposal for the structure in
this case, but we will make
a slightly different proposal in our main conjecture below.

Another motivation for our conjecture is the observation by
Ruan that in his construction (and also in Zharkov's
construction), the discriminant locus
is built out of amoebas, in fact, out of amoebas with the maximum number
of holes.
Since such amoebas arise from moment maps which
are $2$-to-$1$ over the interior and $1$-to-$1$ over the edges, it
is natural to identify the set of singular points of $\pi$
with the algebraic
curve whose moment map image is
 the amoeba.  This is what we do in our main conjecture.

A third input to our conjecture is
the Harvey--Lawson fibration of $\mathbb{C}^3$, which is the standard
model of a ``positive'' vertex with a Lagrangian structure.  In that 
fibration, the set of
singular points consists of the coordinate axes in $\mathbb{C}^3$, which
meet in a ``transverse triple point.''  We
conjecture that this is a general property of positive vertices.

\begin{figure}
\begin{center}
\includegraphics[scale=0.8]{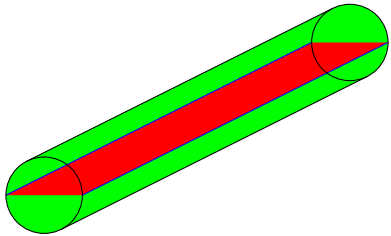}
\end{center}
\caption{Tori which bound a holomorphic disk: the global picture.}
\label{fig:cylinder}
\end{figure}

\begin{conjecture}[Geometry]
\label{con:geometry}
Let $\pi: X\to B$ be a supersymmetric $T^3$
fibration of a compact
Calabi--Yau threefold with respect to a Calabi--Yau metric
whose compatible complex structure is sufficiently close to a
large complex structure limit point,
and whose
K\"ahler class is sufficiently deep in the K\"ahler cone.  Then
$\pi$ is piecewise smooth and
\begin{itemize}
\item[a)] The set $\Sigma\subset X$ of singular points of fibers of $\pi$ is a
complex subvariety of $X$ of complex dimension $1$.

\item[b)] All singular points of $\Sigma$ are transverse triple points, 
locally of the
form $\{z_1z_2=z_1z_3=z_2z_3=0\}$ for local complex coordinates $z_1$,
$z_2$, $z_3$.

\item[c)] For each connected component $\Sigma_\alpha$ of $\Sigma$, 
$\pi(\Sigma_\alpha)$
is contained in a (real) surface $\mathcal A_\alpha\subset B$, and the map
$\pi|_{\Sigma_\alpha}$ is generically $2$-to-$1$ onto its
image $\pi(\Sigma_\alpha)$, which has the topology of
 a compactified amoeba with $g(\Sigma_\alpha)$ holes.

\item[d)] The discriminant locus $\Delta$ retracts to a trivalent
graph $\Gamma$ which
is the union of the spines of the (topological) compactified amoebas
$\pi(\Sigma_\alpha)$.  The graph $\Gamma$ has all of the properties
in Conjecture~\ref{con:topology}.

\item[e)] The positive vertices of the graph $\Gamma$ are the points
in $\Delta$
at which the spines of the various compactified amoebas meet.
The map $\pi$ puts the singular points of $\Sigma$ in one-to-one 
correspondence with the positive vertices

\item[f)] The set $\mathcal H\subset B$ of fibers which are either singular
or contain the boundary of at least one holomorphic disk retracts onto
the discriminant locus $\Delta$.  
(Although in the local example illustrated in
Figure~\ref{fig:thickribbon} the set $\mathcal H$ extended far away from
$\Delta$, we expect that in global examples $\mathcal H$ will be
confined to a small
neighborhood of $\Delta$, as illustrated in Figure~\ref{fig:cylinder}.)
The map $\pi$ is 
smooth\footnote{The referee points out that 
due to the
smoothness of local moduli for special Lagrangian submanifolds 
\cite{MR1664890}, we should even expect $\pi$
to be smooth on the larger set $X-\pi^{-1}(\Delta)$.}
when restricted to
$X-\pi^{-1}(\mathcal{H})$.
\end{itemize}
\end{conjecture}

We have been deliberately vague about the notion of a topological amoeba
and its spine, since we don't know how much of the theory 
of amoebas should be expected
to go through.  It would be very interesting to know, for example, if
some version of the Ronkin function can be defined for a supersymmetric
$T^3$ fibration.

Note that one of the things which could happen if we attempt to
deform this structure too far away from a large complex structure limit
point is that $\pi|_{\Sigma_\alpha}$ might stop being generically
$2$-to-$1$, as happens for moment maps of algebraic curves.  It would
be very interesting to see what happens to supersymmetric $T^3$ fibrations
in that case.  Presumably, something more general than Joyce's phenomenon
of two singular points per fiber is going on here (as Joyce briefly
discusses in \cite[Section 8.2]{Joyce:SYZ}).

\begin{figure}
\begin{center}
\includegraphics[scale=0.7]{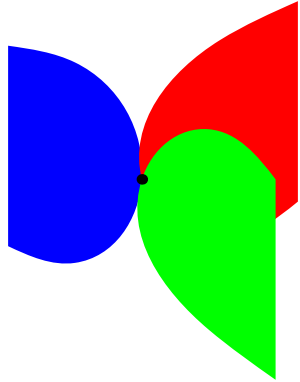}
\end{center}
\caption{The discriminant locus near a positive vertex.}\label{fig:tripledisc}
\end{figure}

Among the consequences of our conjecture is a specific prediction
for the structure of $\Delta$ near a positive vertex.  Using local
affine coordinates to identify a neighborhood of the vertex with the
first cohomology of the fiber, the three local pieces
of $\Delta$ must be contained in the three monodromy-invariant $2$-planes,
which meet along a common line but are distinct.  However, because the
corresponding singular point of $\Sigma$ is a transverse triple point,
the thickening of each piece of 
the discriminant locus will need to ``thin down'' near the positive vertex
so that the three pieces of $\Delta$ meet in a single point,
leading to a description of the discriminant locus similar to that
illustrated in Figure~\ref{fig:tripledisc}.  (This ``thinning down''
was absent from Joyce's proposal about the positive vertices.)
Notice that, as Joyce observed, 
the discriminant locus has a markedly different local structure
near positive and negative vertices and therefore we cannot hope for 
a mirror symmetry statement which simply dualizes all nonsingular tori
in the fibration.

The conjecture that $\pi$ is smooth on a region whose complement retracts
to $\Delta$ was motivated in part by the properties of the construction
of Casta\~no-Bernard and Matessi
\cite{math.SG/0611139}.  
The conjecture that $\mathcal{H}$ provides such a region
is motivated in part by Joyce's observation that---at least
in examples---the boundary between the set of tori bounding holomorphic
disks and the set not bounding holomorphic disks is a boundary along 
which $\pi$ fails to be smooth.  
An additional motivation for part (f) of the conjecture is the hope that
a proper understanding of the disk
contributions to the physical ``T-duality'' construction will restore 
the symmetry between the $T^3$ fibrations on the original Calabi--Yau
manifold and its mirror partner: the fibrations would consist of dual tori on
the complement of $\mathcal{H}$ (where Hitchin's Legendre transform will
relate the affine structures), with the duality between the tori
somehow modified within $\mathcal{H}$ by the
disk contributions.

\section{Degenerations}

We close with a final conjecture, which is perhaps less well-motivated
than Conjecture~\ref{con:geometry}, but which proposes an explanation
for why the
structures we expect from supersymmetric $T^3$ fibrations are 
related to the complex structure
being near a large complex structure limit point (as the
physics suggests).  This final conjecture also points the way towards a
connection between our conjectures and the interesting program
of Gross and Siebert
\cite{GrossSiebert:log,GrossSiebert:logI,GS:affinecomplex,GS:mslogII},
which formulates mirror symmetry in terms of degenerations of
algebraic varieties.

Our final conjecture essentially says that the algebraic curves $\Sigma_\alpha$
should arise from a large complex structure degeneration of the 
Calabi--Yau threefolds.

\begin{conjecture}[Degeneration]
Let $\mathcal X\to S$ be a proper flat family of threefolds whose generic point
$X_\eta$ is a Calabi--Yau threefold, and whose fiber $X_{s_0}$ at
some special point $s_0\in S$ is a large complex structure degeneration of
the form $X_{s_0} = \bigcup X_j$, where the $X_j$ are the components
of $X_{s_0}$.  Equip $\mathcal X\to S$ with a relative
K\"ahler metric $g$ whose K\"ahler classes are sufficiently deep in
 the K\"ahler cone.  
Then there exist (non-flat) families of subvarieties
$\mathcal C_{jk}\subset \mathcal X$ such that $(C_{jk})_{s_0}=X_j\cap X_k$,
but $(C_{jk})_s$ is nonsingular of complex dimension $1$ when $X_s$ is
nonsingular,\footnote{Note that the map $\mathcal C_{jk}\to S$ is not equidimensional: the
fiber over $s_0$ is a surface while the fiber over a general point $s$
is a curve.} 
such that for all $s$ sufficiently close to $s_0$ there is
a supersymmetric $T^3$ fibration of $X_s$ with respect to $g_s$
whose singular locus
$\Sigma_s$ is precisely $\bigcup (C_{jk})_s$.
\end{conjecture}

This is the structure found in the case of Calabi--Yau hypersurfaces
in toric fourfolds: in that case, each $X_j\cap X_k$ is an
intersection of toric divisors,
which meets the nearby nonsingular Calabi--Yau threefolds  $X_s$
in a complex curve $(C_{jk})_s$; the union of
those curves, in the constructions
of Zharkov and of Ruan,
forms the set $\Sigma_s$
of singular points of $X_s$.

\subsection*{Acknowledgments}

It is a pleasure to thank
Emanuel Diaconescu, 
Robbert Dijkgraaf, 
Mark Gross, 
Christian Haase, 
Dominic Joyce, 
Ronen Plesser,
Cumrun Vafa,
and
Ilia Zharkov 
for useful conversations about the topics reported on here,
and to thank the referee for some helpful comments.

I would also like to thank Kansas State University for hosting
the conference on Tropical Geometry and Mirror Symmetry, 
and the Aspen Center for Physics where the
writing of this paper was completed.

\providecommand{\bysame}{\leavevmode\hbox to3em{\hrulefill}\thinspace}

\end{document}